\documentclass[12pt]{article}
\usepackage{amsthm,amssymb,latexsym,amsxtra}

\theoremstyle{plain}
\newtheorem{thm}{Theorem}[section]
\newtheorem{cor}[thm]{Corollary}
\newtheorem{lem}[thm]{Lemma}
\newtheorem{pro}[thm]{Proposition}

\theoremstyle{definition}
\newtheorem{dfn}[thm]{Definition}

\newtheorem{block}[thm]{}

\begin{document}
\title{\bf Hurwitz spaces of Galois coverings of $\mathbb{P}^1$ with Galois groups Weyl groups}
\date{}
 \author{Vassil Kanev}
\maketitle

\def\thefootnote{}
\footnote{Research supported by the italian MIUR, research program
"Geometria sulle variet\`{a} algebriche"}

\begin{abstract}
We prove the irreducibility of the Hurwitz spaces which parametrize 
Galois coverings of $\mathbb{P}^1$ whose Galois group is an arbitrary 
Weyl group and the local monodromies are reflections. This generalizes 
a classical theorem due to Clebsch and Hurwitz.
\end{abstract}

\section*{Introduction}
A classical theorem due to Clebsch and Hurwitz states that the Hurwitz 
space $\mathcal{H}_{d,n}(\mathbb{P}^1)$, which parametrizes 
irreducible coverings of $\mathbb{P}^1$ of degree $d$ simply ramified 
in $n$ points, is irreducible \cite{Hur}. 
Coverings of $\mathbb{P}^1$ of degree $d$ with a fixed monodromy group 
$G\subset S_d$, the related Galois coverings with Galois group $G$, and 
the corresponding Hurwitz spaces were studied in connection with the 
inverse Galois problem (see \cite{Vo} and the references therein). The 
irreducibility of such spaces is a relevant problem for this theory 
and was verified in few cases \cite{BF}, \cite{FB}, \cite{FV}.

Replacing $S_d$ by an arbitrary Weyl group one may ask whether the theorem of
Clebsch and Hurwitz could be generalized to Galois coverings of 
$\mathbb{P}^1$ with Galois groups  Weyl groups. Coverings of this type 
are interesting on their own. They appear in the study of spectral 
curves, integrable systems, generalized Prym varieties, Prym-Tyurin 
varieties \cite{Do, K1, K2, Sc}. The generalized Prym maps yield morphisms from
the Hurwitz spaces of coverings with monodromy groups contained in a Weyl group
to Siegel modular varieties which parametrize Abelian varieties with a fixed 
polarization type. If one proves the irreducibility and the unirationality of some 
Hurwitz spaces and the dominance of the Prym maps, this would imply the unirationality
of the corresponding Siegel modular varieties. This idea was successfully realized
in proving the unirationality of $\mathcal{A}_3(1,1,d)$ and $\mathcal{A}_3(1,d,d)$ for $d\leq 4$
\cite{K3}, \cite{K4} (the case $d=5$ is a work in progress). We hope the method may be
extended considering coverings with monodromy groups an arbitrary irreducible Weyl group.
The irreducibility of the Hurwitz spaces is the first issue to address here. For coverings 
of $\mathbb{P}^1$ Hurwitz showed in \cite{Hur} that the problem of irreducibility is 
reduced to a purely algebraic problem about transitivity of certain actions of the braid groups on tuples 
of elements of the monodromy group (see Section~\ref{s1} for details). We mention that 
analogous reduction, involving braid groups of Riemann surfaces of positive genus, was found 
in \cite{K5} for Hurwitz spaces of coverings of a fixed smooth, projective curve of positive genus.

In the present paper we generalize the result of Clebsch and Hurwitz and prove 
in Theorem~\ref{s2.31} the irreducibility of the Hurwitz spaces which 
parametrize Galois coverings of $\mathbb{P}^1$ whose Galois group is an 
arbitrary Weyl group $W$ and which have simple branching in the sense that 
every local monodromy is a reflection. We notice that when $W$ is the 
Weyl group of an irreducible root system of type $D_{r},B_{r}$ or $C_{r}$ 
the result is already known. The case $W(D_{r})$ is treated in 
\cite{BF}. The case $W(B_{r})$, which is the same as $W(C_{r})$, is easily 
reduced to the theorem of Clebsch and Hurwitz. One consequence of our result is the 
topological classification of the 
coverings we consider. Namely, Clebsch gave a normal form for the local 
monodromies of a simple covering $\pi : X\to \mathbb{P}^1$ (\cite{Cl}, 
cf. \cite{Fu} proof of Proposition~1.5). Our Corollary~\ref{s2.32} 
gives a normal form for the local monodromies when the monodromy group 
is a Weyl group and the branching is simple.

We mention two other recent papers where the problem of irreducibility of Hurwitz spaces
of coverings of $\mathbb{P}^1$ was studied. Let $W$ be a finite irreducible Coxeter group of rank $r$.
S. Humphries studied in \cite{Hum} Hurwitz actions of the $r$-strand braid group on $r$-tuples of 
reflections of W and 
one of his results has the following corollary.
Let $W$ be of type $A_r,B_r,D_r,E_6,E_7,F_4, I_3$ or $I_4$. 
Theorem~1.2 (ibid.) implies the irreducibility of the Hurwitz space parametrizing
the irreducible Galois covers of $\mathbb{P}^1$ with Galois group $W$, branched in $r+1$ points, 
in $r$ of which the local monodromies are reflections. F. Vetro studied in \cite{Ve} coverings of 
$\mathbb{P}^1$ of degree $2r$ whose monodromy group is contained in $W(B_r)\subset S_{2r}$. She 
proved the irreducibility of the corresponding Hurwitz space when the local monodromies at all branch
points, except possibly one, are reflections. We do not know whether our theorems 2.5 and 2.7 remain 
valid if one replaces Weyl groups by finite Coxeter groups. While only Weyl groups are relevant 
for the Siegel modular varieties, a possible generalization to finite Coxeter groups might be of interest
for the inverse Galois theory (see e.g. \cite{Vo}).

\section{Hurwitz spaces and Weyl groups}\label{s1}

In \S\ref{s1.5}--\S\ref{s1.9ac} we recall some facts about Hurwitz 
spaces. The references for this material are \cite{Fu}, \cite{Fr}, 
 \cite{Vo}.

\begin{block}\label{s1.5}
Let $\pi :X\to \mathbb{P}^{1}$ be a Galois covering with Galois  group 
$G$. We assume $G$ acts on $X$ on the left. We call $\pi$ a 
$G$-covering for short. Let $D\subset \mathbb{P}^1$ be the 
discriminant locus of $\pi$ and let $b_0\in \mathbb{P}^1-D$. We 
consider the fundamental groups with 
multiplication defined by means of composition of arcs 
$\gamma_{1}*\gamma_{2}$  where one first travels along $\gamma_{1}$ 
and then along $\gamma_{2}$. Let $x_0\in \pi^{-1}(b_0)$. The monodromy 
homomorphism $m_{x_0}:\pi_1(\mathbb{P}^1-D,b_0)\to G$ is defined as 
follows. If $\gamma$ is a closed arc in $\mathbb{P}^1-D$ based at 
$b_0$, let $\hat{\gamma}:[0,1]\to X-\pi^{-1}(D)$  be its lifting which 
starts at $x_0$ and ends at $gx_0$. Then $m_{x_0}([\gamma])=g$. We 
notice that if instead one considers Galois coverings where $G$ acts 
on the right, then $m_{x_0}([\gamma])=g$ if $\hat{\gamma}(0)=x_0, 
\hat{\gamma}(1)=x_0g^{-1}$.

Let $\ell$ be an arc which connects $b_0$ with a point $b\in D$ and 
contains none of the other points of $D$. Let $\gamma$ be a closed arc 
which begins at $b_0$, travels along $\ell$ to a point near $b$, makes 
a small counterclockwise loop around $b$ and returns to $b_0$ along 
$\ell$. The element $t=m_{x_0}([\gamma])$ is the \emph{local 
monodromy} around $b$ along $\gamma$. The conjugacy class 
$\{gtg^{-1}|g\in G\}$ does depend neither of the choice of $\gamma$, nor 
of the choice of $x_0\in \pi^{-1}(b_0)$, nor of the choice of $b_0\in 
\mathbb{P}^1-D$. It characterizes the ramification type at the 
discriminant point $b\in D$. 

Let $n=|D|$. An \emph{arc system} is a collection of $n$ simple arcs 
($=$ embedded intervals) which join $b_0$ with the points of $D$ and 
do not meet outside $b_0$. One defines an ordering of an arc system by 
choosing arbitrarily the first one and numbering the arcs by the 
directions of departure in counterclockwise order. One associates to 
an arc system $n$ closed arcs as above. We call the obtained 
$\gamma_{1},\ldots ,\gamma_{n}$ a \emph{standard system} of closed 
arcs. Their homotopy classes generate $\pi_1(\mathbb{P}^1-D,b_0)$ with 
the only relation $\gamma_{1}*\cdots *\gamma_{n}\simeq 1$. 
Let $t_i=m_{x_0}([\gamma_{i}])$. One has $t_{1}\cdots t_{n}=1$. 
Conversely, let the $n$-tuple $(t_{1},\ldots,t_{n})$ of elements of 
$G$ satisfies $t_i\neq 1$ for $\forall i$ and $t_{1}\cdots t_{n}=1$. Let 
$D\subset \mathbb{P}^1-b_0$ be an arbitrary set of $n$ points, let 
$\gamma_{1},\ldots ,\gamma_{n}$ be a simple system of closed arcs 
constructed as above and let $m:\pi_1(\mathbb{P}^1-D,b_0)\in G$ be the 
homomorphism defined by $m_{x_0}([\gamma_{i}])=t_i$. Then by Riemann's 
existence theorem there is a $G$-covering $\pi :X\to \mathbb{P}^{1}$ 
branched in $D$, and a point $x_0\in \pi^{-1}(b_0)$ such that 
$m_{x_0}=m$. The topological cover $X-\pi^{-1}(D)$ is connected, 
equivalently $X$ is irreducible, if and only if $t_1,\ldots,t_n$ 
generate $G$. 
\end{block}
\begin{dfn}\label{s1.7}
An $n$-tuple $(t_1,\ldots,t_n)$ of elements of a group $G$ which 
satisfy $t_i\neq 1$ for $\forall i$ and $t_{1}\cdots t_{n} = 1$ is 
called a \emph{Hurwitz system}.
\end{dfn}
\begin{block}\label{s1.7a}
Two $G$-coverings $\pi :X\to \mathbb{P}^{1}$ and $\pi' :X'\to 
\mathbb{P}^{1}$ are called equivalent if there is a $G$-equivariant 
isomorphism $f:X\to X'$ such that $\pi'=\pi\circ f$. Suppose 
furthermore $\pi$ and $\pi'$ are not branched in $b_0\in 
\mathbb{P}^1$. Let $x_0\in \pi^{-1}(b_0), x'_0\in \pi^{-1}(b_0)$. 
The pairs $(\pi :X\to \mathbb{P}^{1},x_0)$ and $(\pi' :X'\to 
\mathbb{P}^{1},x'_0)$ are called equivalent if there is an isomorphism $f$ as 
above which satisfies furthermore $f(x_0)=x'_0$.

Let $C_1,\ldots,C_k$ be conjugacy classes of $G$,  $C_i\neq C_j$ 
if $i\neq j$. Let $\underline{n}=n_1C_1+\cdots +n_kC_k$  be a formal 
sum where $n_i\in \mathbb{N}$. Let $|\underline{n}|=n_1+\cdots +n_k$. 
In this paper we study two types of Hurwitz spaces which we first 
define as sets. The points of 
$\mathcal{H}_{G;\underline{n}}(\mathbb{P}^1)$ are the equivalence 
classes of $G$-coverings $[\pi: X\to \mathbb{P}^1]$ with $n=|\underline{n}|$ 
discriminant points such that $n_i$ of these points have local 
monodromies belonging to $C_i, i=1,\ldots,k$, and moreover $X$ is irreducible. 
The points of $\mathcal{H}_{G;\underline{n}}(\mathbb{P}^1,b_0)$ are 
the equivalence classes of pairs $[\pi: X\to \mathbb{P}^1,x_0]$ where 
$\pi$ is as above and furthermore it is unramified at $b_0\in 
\mathbb{P}^1$ and $\pi(x_0)=b_0$. The following properties are known.

(i)\quad $\mathcal{H}_{G;\underline{n}}(\mathbb{P}^1)\neq \emptyset 
\Leftrightarrow \mathcal{H}_{G;\underline{n}}(\mathbb{P}^1,b_0)\neq \emptyset 
\Leftrightarrow $ there exists a Hurwitz system $(t_1,\ldots,t_n)\in G^n,\; 
t_{1}\cdots t_{n}=1$, such that $t_1,\ldots,t_n$ generate $G$ (cf. 
\S\ref{s1.5}). We assume for the next properties that the Hurwitz 
spaces are nonempty.

(ii)\quad $\mathcal{H}_{G;\underline{n}}(\mathbb{P}^1)$ has a 
canonical complex analytic structure such that the map given by $[X\to 
\mathbb{P}^1]\mapsto D$ is a finite, \`{e}tale holomorphic map 
$\mathcal{H}_{G;\underline{n}}(\mathbb{P}^1)\to 
(\mathbb{P}^1)^{(n)}-\Delta$. Here $\Delta$ is the codimension one 
subvariety consisting of effective non-simple divisors of degree $n$. 
Similarly $\mathcal{H}_{G;\underline{n}}(\mathbb{P}^1,b_0)$ has a 
canonical complex analytic structure and a finite \`{e}tale covering 
$\mathcal{H}_{G;\underline{n}}(\mathbb{P}^1,b_0)\to 
(\mathbb{P}^1-b_0)^{(n)}-\Delta$.

(iii)\quad The complex analytic spaces and the coverings of (ii) are 
algebraic.

(iv)\quad Let us fix a $D\in (\mathbb{P}^1-b_0)^{(n)}$ and let us 
choose a standard system of closed arcs $\gamma_{1},\ldots 
,\gamma_{n}$. Varying $m:\pi_1(\mathbb{P}^1-D,b_0)\to G$ let us identify 
the fiber of $\mathcal{H}_{G;\underline{n}}(\mathbb{P}^1,b_0)\to 
\linebreak (\mathbb{P}^1-b_0)^{(n)}-\Delta$ over $D$ with the set of Hurwitz 
systems $(t_1,\ldots,t_n)\in G^n,\; t_{1}\cdots t_{n}=1$, such that $n_i$ 
of its elements belong to $C_i,\; i=1,\ldots,k$, and 
$t_1,\ldots,t_n$ generate $G$. Consider the monodromy action of the 
braid group $\pi_1((\mathbb{P}^1-b_0)^{(n)}-\Delta,D)$ on the fiber 
over $D=\{b_1,\ldots,b_n\}$. Then the action of the elementary braids 
$\sigma_{i},{\sigma_{i}}^{-1}$ which fix $b_j$ for $j\neq i,i+1$ and 
rotate $b_i,b_{i+1}$ at angles $\pm\pi$, is given by the following 
formulae

\begin{equation}\label{es1.8}
\begin{split}
&(t_1,\ldots,t_{i-1},t_{i},t_{i+1},t_{i+2},\ldots,t_n)
\mapsto
(t_1,\ldots,t_{i-1},t_{i}t_{i+1}t_{i}^{-1},t_{i},t_{i+2},\ldots,t_n)\\
&
(t_1,\ldots,t_{i-1},t_{i},t_{i+1},t_{i+2},\ldots,t_n)
\mapsto
(t_1,\ldots,t_{i-1},t_{i+1},t_{i+1}^{-1}t_{i}t_{i+1},t_{i+2},\ldots
,t_n)\\
\end{split}
\end{equation}
We call such transformations of $n$-tuples \emph{elementary 
transformations} or \emph{braid moves}. They determine uniquely the 
monodromy action of the braid group since the elementary braids 
generate it.

(v)\quad The image of the forgetful map $[X\to 
\mathbb{P}^1,x_0]\mapsto [X\to \mathbb{P}]$ is a Zariski open, dense 
subset $\mathcal{U}(b_0)\subset 
\mathcal{H}_{G;\underline{n}}(\mathbb{P}^1)$ consisting of $[X\to 
\mathbb{P}]$ which are unramified at $b_0$. Let us denote by 
$[t_1,\ldots ,t_n]$ the orbit of an $n$-tuple of elements of $G$ with 
respect to the conjugacy action $(t_1,\ldots,t_n)\mapsto 
(st_1s^{-1},\ldots,st_ns^{-1}),\; s\in G$. Then the fiber of 
$\mathcal{U}(b_0)\to (\mathbb{P}^1-b_0)^{(n)}-\Delta$ over $D$ may be 
identified with the set $\{[t_1,\ldots ,t_n]\}$ where 
$(t_1,\ldots,t_n)$ runs over all Hurwitz systems satisfying the 
conditions of (iv). This set is called \emph{Nielsen class} and 
denoted by $Ni(\underline{n},G)$. The monodromy action of the braid 
group $\pi_1((\mathbb{P}-b_0)^{(n)}-\Delta,D)$ on the Nielsen class is 
determined by formulae \eqref{es1.8}.
\end{block}

\begin{dfn}\label{s1.9a}
We call two $n$-tuples of elements of $G$ braid-equivalent if one can 
be obtained from the other by a finite sequence of elementary 
transfromations \eqref{es1.8}. We denote the braid equivalence by $\sim$. 
\end{dfn}

\begin{block}\label{s1.9ac}
Using (ii), (iv) and (v) of \S\ref{s1.7a} one obtains the following 
conclusion.

\emph{The Hurwitz space 
$\mathcal{H}_{G;\underline{n}}(\mathbb{P}^1,b_0)$ is irreducible if 
and only if every two Hurwitz systems satisfying the conditions of 
(iv) are braid-equivalent. The Hurwits space 
$\mathcal{H}_{G;\underline{n}}(\mathbb{P}^1)$ is irreducible if and 
only if every two $G$-orbits of Hurwitz systems satisfying the 
conditions of (v) are braid-equivalent}
\end{block}

\begin{lem}[\cite{FB} p.102, \cite{Vo} Lemma 9.4]\label{s1.9aa}
Let $(t_1,\ldots,t_n)\in G^n,\; t_{1}\cdots t_{n}=1$ be a Hurwitz 
system. Then
\renewcommand{\theenumi}{\roman{enumi}}
\begin{enumerate}
\item $(t_1,t_2,\ldots,t_n)\sim (t_2,\ldots,t_n,t_1)$.
\item If $s\in \langle t_{1},\ldots ,t_{n} \rangle \subset G$ then
$(t_1,\ldots,t_n)\sim (st_1s^{-1},\ldots,st_ns^{-1})$
\end{enumerate}
\end{lem}

\begin{cor}\label{s1.9ad}
The forgetful map $\mathcal{H}_{G;\underline{n}}(\mathbb{P}^1,b_0)\to 
\mathcal{H}_{G;\underline{n}}(\mathbb{P}^1)$ establishes a bijective 
correspondence between the connected (= irreducible) components of the 
two Hurwitz spaces
\end{cor}

\begin{lem}\label{s1.9ab}
Suppose the $n$-tuple $(t_1,\ldots,t_n)\in G^n$ contains the adjacent 
pair $t_k=t,t_{k+1}=t^{-1}$. Then performing a sequence of elementary 
transformations \eqref{es1.8} one may move the pair $(t,t^{-1})$ at 
any two consequitive places without changing the other elements of the 
$n$-tuple.
\end{lem}
\begin{proof}
This follows from the braid equivalences $(u,t,t^{-1})\sim 
(t,t^{-1}ut,t^{-1})\sim (t,t^{-1},u)$ and $(t,t^{-1},u)\sim 
(t,t^{-1}ut,t^{-1})\sim (u,t,t^{-1})$.
\end{proof}

\begin{lem}\label{s1.9b}
Let $(t_1,\ldots,t_n)$ be an $n$-tuple of elements of $G$ such that 
$t_{i+1}=t_{i}^{-1}$. Let $H$ be the subgroup generated by 
$t_1,\ldots,t_{i-1},t_{i+2},\ldots,t_n$. Then for every $h\in H$ one has 
\[
(t_1,\ldots,t_{i-1},ht_{i}h^{-1},ht_{i+1}h^{-1},t_{i+2},\ldots,t_n)
\sim
(t_1,\ldots,t_n).
\]
\end{lem}
\begin{proof}
Let $H_1$ be the subset consisting of elements $h$ such that 
\[
(t_1,\ldots,t_{i-1},hth^{-1},ht^{-1}h^{-1},t_{i+2},\ldots,t_n)
\sim
(t_1,\ldots,t_{i-1},t,t^{-1},t_{i+2},\ldots,t_n)
\]
holds for every $t\in G$. By reflexivity and transitivity of $\sim$ it 
follows that $H_1$ is a subgroup of $G$. The statement of the lemma 
would be proved if we show that $t_{j}\in H_{1}$ for every $j\neq 
i,i+1$. Using Lemma~\ref{s1.9ab} we move $(t,t^{-1})$ to the right of 
$t_{j}$. Then we have $(t_{j},t,t^{-1})\sim 
(t_jtt_{j}^{-1},t_{j},t^{-1})\sim 
(t_jtt_{j}^{-1},t_jt^{-1}t_{j}^{-1},t_{j})$. We then move the pair $(
t_jtt_{j}^{-1},t_jt^{-1}t_{j}^{-1})$ back to the initial position.
\end{proof}

\begin{block}\label{s1.9c}
Let $R$ be a root system in a real vector space $V$ (see \cite{Bo} Ch. VI 
or \cite{Hu1} Ch.III). Let $r=\dim V$ be the rank of $R$. We assume $R$ 
is reduced, i.e. for each $\alpha \in R$ one has $R\alpha \cap 
R=\{\alpha,-\alpha\}$. Let $W=W(R)$ be the Weyl group generated by the 
reflections $\{s_{\alpha}|\alpha \in R$\}. Let $(\; |\; )$ be a 
$W$-inveriant inner product of $V$. Following \cite{Bo} we denote 
$\alpha^{\vee}=\frac{2\alpha}{(\alpha|\alpha)}$ and 
$n(x,\alpha)=\frac{2(x|\alpha)}{(\alpha|\alpha)}$. Then 
$s_{\alpha}=x-n(x,\alpha)\alpha$. The values of $n(\alpha,\beta)$ for 
$\alpha,\beta\in R$ are given in \cite{Bo} Ch.VI~\S1.3. Let $C$ be a 
chamber, $B(C)=\{\alpha_{1},\ldots,\alpha_{r}\}$ be the corresponding 
base of $R$, $R=R^{+}\cup R^{-}$ be the decomposition into positive 
and negative roots. Every total ordering of $V$ determines uniquely a 
base of $R$ composed of the set of positive roots indecomposable into 
sums of other positive roots \cite{Hu2} \S1.3 and \cite{Bo} Ch.VI~\S7. 
Vice versa given a base of $R$ and choosing its linear ordering one may 
consider the corresponding lexicographic ordering of $V$. The simple 
system associated with this total ordering is the given base of $R$. 
Suppose $R$ is an irreducible root system. If $R$ is simply laced, 
i.e. of type $A_{r}, D_{r}, E_{6}, E_{7}$ or $E_{8}$ then it consists 
of one $W$-orbit. If $R$ is non-simply laced, i.e. of type $B_{r}, 
C_{r}, F_{4}$ or $G_{2}$, then it consists of two $W$-orbits 
$R=R_{s}\cup R_{_{\ell}}$ called short and long roots respectively. 
\end{block}

\begin{dfn}\label{s1.9d}
Let $R$ be a root system and let $W$ be its Weyl group. A Galois 
covering $\pi : X\to \mathbb{P}^1$ with Galois group $W$ is called 
simply ramified if every local monodromy is a reflection.
\end{dfn}

\begin{block}\label{s1.9da}
Simply ramified $W$-coverings of $\mathbb{P}^1$ yield Hurwitz systems 
$t_{1}\cdots t_{n}=1$ where $t_{i}$ are reflections in $W$. Applying 
the canonical homomorphism $\epsilon:W\to \{1,-1\}$ we see that $n$ is 
even for these Hurwitz systems.

Let $R=R^{(1)}\sqcup \ldots \sqcup R^{(k)}$ be the decomposition into 
disjoint union of irreducible root systems. Let $W^{(i)}=W(R^{(i)})$ 
be the corresponding Weyl groups. Every reflection in $W$ belongs to 
some $W^{(i)}$ and every conjugacy class of reflections in $W$ is a 
conjugacy class of reflections in some $W^{(i)}$ with respect to 
$W^{(i)}$. Simplifying the notation of \S\ref{s1.7a} we may specify 
the branching data of a simply ramified $W$-covering $\pi : X\to 
\mathbb{P}^1$ by 
$\underline{n}=(\underline{n}^{(1)},\ldots,\underline{n}^{k})$ where: 
if $R^{(i)}$ is simply laced $\underline{n}^{(i)}=n^{(i)}$ and denotes 
the number of discriminant points whose local monodromies are 
reflections in $W^{(i)}$; if $R^{(i)}$ is non-simply laced then 
$\underline{n}^{(i)}=(n_{s}^{(i)},n_{_{\ell}}^{(i)})$ where 
$n_{s}^{(i)}$, respectively $n_{_{\ell}}^{(i)}$, denotes the number of 
discriminant points whose local monodromies are reflections with respect to 
short roots, respectively long roots, in $W^{(i)}$. Our aim is to 
prove that the Hurwitz spaces 
$\mathcal{H}_{W;\underline{n}}(\mathbb{P}^1,b_0)$ and 
$\mathcal{H}_{W;\underline{n}}(\mathbb{P}^1)$ are irreducible by means 
of studying braid equivalences between Hurwitz systems of reflections 
in $W$.
\end{block}

\begin{dfn}\label{s1.10}
Let $\{t_{1},\ldots,t_{m}\}$ be a subset of a group $G$. Replacing a 
pair $t_{i},t_{j}$ by $t_{i},t_{i}t_{j}t_{i}^{-1}$ one obtains a new 
set $\{t'_{1},\ldots,t'_{m}\}$ where $t'_{k}=t_{k}$ if $k\neq j$ and 
$t_{j}=t_{i}t_{j}t_{i}^{-1}$. This transformation of subsets of a group is 
called 
\emph{Nielsen transformation}. 
\end{dfn}
The following theorem is a particular case of a well-known result 
\cite{Co,De,Dy}. The part concerning Nielsen transformations, which 
we need, is treated only in \cite{Be} in a special case, so we include 
a simple proof.

\begin{thm}\label{s1.11}
Let $R\subset V$ be a root system, let $W=W(R)$ be its Weyl group, let 
$\geq$ be a total ordering of $V$ and let 
$\{\alpha_{1},\ldots,\alpha_{r}\}$ be the corresponding base of $R$. 
Let $T=\{t_{1},\ldots,t_{m}\}\subset W$ be a set of reflections, 
$t_{i}=t_{\beta_{i}}$ with $\beta_{i} \in R^{+}$. Let $W'$ be the subgroup 
generated by $T$, let $R'=W'\cdot \{\beta_{1},\ldots,\beta_{m}\}$ and 
let $V'$ be the span of $R'$. Then $R'$ is a root system in $V'$ and 
$W'$ is its Weyl group. Furthermore if 
$\{\alpha'_{1},\ldots,\alpha'_{_{\ell}}\}$ is the base of $R'$ associated 
with the total ordering induced by $\geq$ on $V'$, then the set 
$\{s_{\alpha'_{1}},\ldots,s_{\alpha'_{_{\ell}}}\}$ may be obtained from 
$T$ by a finite sequence of Nielsen transformations.
\end{thm}

\begin{proof}
First we notice that if $s_{\alpha},s_{\beta}$ are reflections in $W$ 
then $s_{\alpha}s_{\beta}s_{\alpha}=s_{s_{\alpha}(\beta)}$. Let $(\; 
|\; )$ be a $W$-invariant inner product of $V$. For each $\alpha \in 
R,\; \alpha = a_{1}\alpha_{1}+\cdots +a_{r}\alpha_{r}$, let 
$ht(\alpha)=a_{1}+\cdots +a_{r}$. If $s=s_{\alpha}$ is a reflection 
with $\alpha \in R^{+}$ let $h(s)=ht(\alpha)$. We extend this 
definition to sets of reflections $A=\{\sigma_{i}\}$ letting 
$h(A)=\sum h(\sigma_{i})$. 

CLAIM. \emph{Suppose $\alpha,\beta \in R^{+}$} and $(\alpha|\beta)>0$. 
then there is a Nielsen transformation 
$\{s_{\alpha},s_{\beta}\}\mapsto \{s_{\alpha'},s_{\beta'}\}$ such that 
$h(\{s_{\alpha},s_{\beta}\})>h(\{s_{\alpha'},s_{\beta'}\})$.

For the proof of the claim we may obviously suppose $R$ is 
irreducible. We may furthermore suppose $\| \beta \|\geq\|\alpha\|$. 
First let $\|\alpha\|=\|\beta\|$. We have 
$s_{\alpha}(\beta)=\beta-\alpha\in R$. If $\beta-\alpha\in R^{+}$ then 
$h(s_{s_{\alpha}(\beta)})=ht(\beta-\alpha)<ht(\beta)$ and we may 
decrease $h$ considering the Nielsen transformation 
$\{s_{\alpha},s_{\beta}\}\mapsto 
\{s_{\alpha},s_{\alpha}s_{\beta}s_{\alpha}\}$. If $\beta-\alpha\in 
R^{-}$, then $s_{\beta}(\alpha)=\alpha-\beta\in R^{+}$ and we consider 
$\{s_{\beta}s_{\alpha}s_{\beta},s_{\beta}\}$. If 
$\|\beta\|>\|\alpha\|$ then $s_{\beta}(\alpha)=\alpha-\beta$ and 
$s_{\alpha}(\beta)=\beta-c\alpha$ where $c=2$ if $R$ is of type 
$B_{n},C_{n}$ or $F_{4}$ and $c=3$ if $R$ is of type $G_{2}$. If 
$\alpha-\beta\in R^{+}$ we consider as above 
$\{s_{\beta}s_{\alpha}s_{\beta},s_{\beta}\}$. If 
$s_{\alpha}(\beta)=\beta-c\alpha \in R^{+}$ then we consider 
$\{s_{\alpha},s_{\alpha}s_{\beta}s_{\alpha}\}$. It remains to deal with 
the cases when $\|\beta\|>\|\alpha\|,\; \alpha-\beta \in R^{-},\; 
\beta-c\alpha \in R^{-}$. Suppose first $R$ is of type 
$B_{n},C_{n}$ or $F_{4}$, so $c=2$. Then $\alpha-\beta\in R^{-}$ 
implies $ht(\alpha)<\beta$, so the positive root $2\alpha-\beta$ 
satisfies $ht(2\alpha-\beta)<ht(\beta)$.  Therefore 
$h(s_{\alpha}s_{\beta}s_{\alpha})=h(s_{2\alpha-\beta})<h(s_{\beta})$ 
and we may decrease $h$ by the Nielsen transformation 
$\{s_{\alpha},s_{\beta}\}\mapsto 
\{s_{\alpha},s_{\alpha}s_{\beta}s_{\alpha}\}$. Finally, let $R$ be of 
type $G_{2}$. With the assumptions above we have
\begin{align*}
&h(s_{\alpha})+h(s_{\alpha}s_{\beta}s_{\alpha}) = 
ht(\alpha)+ht(3\alpha-\beta) = 4ht(\alpha) - ht(\beta)\\
&h(s_{\beta}s_{\alpha}s_{\beta})+h(s_{\beta}) = ht(\beta-\alpha) + 
ht(\beta) = 2ht(\beta) - ht(\alpha).
\end{align*}
It is impossible that both numbers $\geq ht(\alpha)+ht(\beta)$. This 
proves the claim.

The claim shows that if among the reflections of $T$ there are two  
$t_{i}=s_{\beta_{i}},\: t_{j}=s_{\beta_{j}}$ with 
$\beta_{i},\beta_{j}\in R^{+}$ and $(\beta_{i}|\beta_{j})>0$, then 
performing a Nielsen transformation of $T$ we may decrease $h(T)$. 
Since $h$ assumes a finite number of positive values on sets of $\leq 
m$ elements we conclude that after a finite number of Nielsen 
transformations we obtain a set of reflections 
$T'=\{s_{\alpha'_{1}},\ldots,s_{\alpha'_{_{\ell}}}\}$ where 
$\alpha'_{i}\in R^{+}$ and $(\alpha'_{i}|\alpha'_{j})\leq 0$ for 
$i\neq j$. We have furthermore that $n(\alpha'_{i},\alpha'_{j})\cdot 
n(\alpha'_{j},\alpha'_{i})\in \{0,1,2,3\}$ if $i\neq j$ since 
$\alpha'_{i},\alpha'_{j}\in R$. By the classification of Coxeter 
graphs and Dynkin diagrams (see e.g. \cite{Hu1}) it follows that 
$\alpha'_{1},\ldots,\alpha'_{_{\ell}}$ is a base of a root system. 
Replacing $T$ by a Nielsen transformation does change neither 
$W'=\langle T\rangle$ nor $R'=W'\cdot \{\beta_{1},\ldots,\beta_{m}\}$. 
We conclude that $R'$ is a root system with base $\alpha'_{1},\ldots,\alpha'_{_{\ell}}$ and 
moreover since all $\alpha'_{i}$ are positive roots in $R$ this is the unique 
base of $R'$ associated with the total ordering $\geq$. 
\end{proof}
\begin{cor}\label{s1.13a}
Let $R\subset V$ be a root system, let $W$ be its Weyl group and let 
$S=\{s_{\alpha_{1}},\ldots,s_{\alpha_{r}}\}$ be a set of simple reflections 
corresponding to a base of $R$. Let $T=\{t_{1},\ldots,t_{m}\}$ be a 
set of reflections which generate $W$. Then one can obtain $S$ from 
$T$ by a finite sequence of Nielsen transformations.
\end{cor}
\begin{proof}
Let us choose a linear ordering of the set 
$\{\alpha_{1},\ldots,\alpha_{r}\}$ and let $\geq$ be the corresponding 
lexicographic ordering of $V$. Using the notation of the theorem we 
have $W'=W,\; R'=R$ (see \cite{Hu2} Ex. 1.14), $V'=V$ and 
$\{\alpha'_{1},\ldots,\alpha'_{_{\ell}}\}=\{\alpha_{1},\ldots,\alpha_{r}\}$ 
since the total ordering $\geq$ determines uniquely a base of 
$R=R'$.
\end{proof}
\section{Irreducibility of Hurwitz spaces}\label{s2}
\begin{block}\label{s2.14}
Suppose $R\subset V$ is an irreducible root system. Let us normalize 
the $W$-invariant inner product $(\; |\; )$ so that 
$(\alpha|\alpha)=2$ for every root $\alpha$ if $R$ is simply laced and 
$(\alpha|\alpha)=2$ for every short root if $R$ is non-simply laced. 
In the latter case for every long root $\beta$ one has 
$(\beta|\beta)=4$ if $R$ is of type $B_{r},C_{r}$ or $F_{4}$ and 
$(\beta|\beta)=6$  if $R$ is of type $G_{2}$. Let us choose a chamber 
$C\subset V$. Let $\lambda$ be the dominant root if $R$ is simply 
laced and let $\lambda$ be the dominant short root if $R$ is 
non-simply laced. It is known that 
\emph{if $\alpha\in R^{+}$ and $\alpha\neq \lambda$ then 
$(\lambda|\alpha^{\vee})$ equals 0 or 1}
(cf. \cite{Bo} Ch.VI~\S1.3).  With the fixed normalization of 
$(~|~)$ we have that $\alpha^{\vee}=\alpha$ if $R$ is of type 
$A_r,D_r,E_6,E_7,E_8$ or if $R$ is of type $B_r,C_r,F_4,G_2$ 
 and $\alpha$ is a short root. So in these cases we have that 
$\alpha\in R^{+}$ and $\alpha\neq \lambda$ implies that 
$(\lambda|\alpha)$ equals 0 or 1. If $R$ is of type $B_r,C_r,F_4$ or 
$G_2$ and $\beta$ is a positive long root then $(\lambda|\beta)$ 
equals 0 or $(\beta|\beta)/2$. 
\end{block}
\begin{pro}\label{s2.15}
Let $R$ be a root system and let $W$ be its Weyl group. Let 
$t_{1}\cdots t_{n}=1$ be a Hurwitz system where $t_{i}$ are 
reflections in $W$. Then $t_{1}\cdots t_{n}=1$ is braid-equivalent to 
a Hurwitz system $t'_{1}t'_{2}\cdots t'_{n}=1$ such that 
$t'_{1}=t'_{2}$.
\end{pro}
\begin{proof}
By Theorem~\ref{s1.11} we may suppose without loss of generality that 
the set $\{t_1,\ldots,t_n\}$ generates $W$. Let 
$R=R^{(1)}\sqcup \ldots \sqcup R^{(k)}$ be the decomposition into 
disjoint union of irreducible root systems. Let $W^{(i)}=W(R^{(i)})$. 
If $\alpha\in R^{(i)},\; \beta\in R^{(j)}$ and $i\neq j$ then 
$(\alpha|\beta)=0$ and $s_{\alpha}$ commutes with $s_{\beta}$. Thus 
one may perform several braid moves to the Hurwitz system $t_{1}\cdots 
t_{n}=1$ to the effect of obtaining a concatenation 
$T^{(1)}T^{(2)}\cdots T^{(k)}=1$ where 
$T^{(i)}$ contains all reflections in $T$ which belong to $W^{(i)}$ in 
the order they appear in $T$. Since $W$ is a direct product of 
$W^{(i)}$ the product of the 
reflections in $T^{(i)}$ equals 1 for each $i$. It suffices to prove 
the proposition for $T^{(1)}$, so we may suppose without loss of 
generality that $R$ is an irreducible root system. The case 
$rk(R)=1,\; W=S_2$ is trivial, so we may furthermore suppose 
$r=rk(R)\geq 2$. If $T=(t_1,\ldots,t_n)$ is braid-equivalent to a 
sequence with two equal reflections then we may move them by 
elementary transformations to the first two places and obtain the 
required Hurwitz system. Let us suppose by way of contradiction that 

\smallskip
(*)\emph{
No sequence braid-equivalent to $(t_1,\ldots,t_n)$ contains two equal 
reflections.
}

\smallskip

Step 1.\; Let $t_{i}=s_{\gamma_{i}}$ with $\gamma_{i}\in R^{+}$. Let 
$\lambda$ be the dominant root defined in \S\ref{s2.14}. If among 
$t_1,\ldots,t_n$ the reflection $s_{\lambda}$ is present we move it to 
the first place by a sequence of elementary transformations of the 
type $(\sigma,\tau)\mapsto (\tau,\tau\sigma\tau)$. Similarly we may 
move to the front all reflections $s_{\gamma_{i}}$ with 
$(\lambda|\gamma_{i})>0$. We obtain a Hurwitz system braid-equivalent 
to $t_{1}\cdots t_{n}=1$ of the form
\begin{equation}\label{es2.16}
s_{\beta_{1}}s_{\beta_{2}}\cdots s_{\beta_{k}}s_{\beta_{k+1}}\cdots 
s_{\beta_{n}} = 1
\end{equation}
where:

(i)\quad $\beta_{i}\in R^{+}$ and $\beta_{i}\neq \beta_{j}$ for $i\neq 
j$ (Assumption~(*));

(ii)\quad $(\lambda|\beta_{i})>0$ for $i=1,\ldots,k$ and 
$(\lambda|\beta_{j})=0$ for $j\geq k+1$;

(iii)\quad if $\lambda\in \{\beta_{1},\ldots,\beta_{k}\}$ then 
$\lambda=\beta_{1}$.

\noindent
Notice that $k\geq 1$ since $W$ cannot be generated by reflections 
$s_{\alpha}$ with $(\lambda|\alpha)=0$. Among the Hurwitz systems 
braid-equivalent to $t_{1}\cdots t_{n}=1$ and satisfying conditions 
(i) -- (iii) we consider those for which 

(iv)\quad $k$ is minimal possible.

\smallskip

Step 2.\; We claim there is a Hurwitz system braid-equivalent to 
$t_{1}\cdots t_{n}=1$ which satisfies conditions (i) -- (iv) of Step 1 
and furthermore 

(v)\quad $(\beta_{i}|\beta_{j})\geq 0$ for every $i,j$ with $1\leq 
i,j\leq k$

\noindent
holds. First suppose $\beta_{1}=\lambda$. Then 
$(\beta_{1},\beta_{i})\geq 0$ for $i\geq 2$ since $\lambda$ is a 
dominant root. Suppose there is a pair $i,j$ with $2\leq i<j\leq k$ 
such that $(\beta_{i}|\beta_{j})<0$. If $\|\beta_{i}\|=\|\beta_{j}\|$ 
then $s_{\beta_{i}}(\beta_{j})=\beta_{j}+\beta_{i}$ and 
$\|s_{\beta_{i}}(\beta_{j})\|=\|\beta_{j}\|$. We have $(\lambda|
s_{\beta_{i}}(\beta_{j}))=(\lambda|\beta_{i})+(\lambda|\beta_{j}) > 
\max\{(\lambda|\beta_{i}),(\lambda|\beta_{j})\}$. From \S\ref{s2.14} 
this is possible only if $\|\beta_{i}\|=\|\beta_{j}\|=2$, and 
$\lambda=\beta_{i}+\beta_{j}$. Performing several braid moves among 
$s_{\beta_{2}},\ldots,s_{\beta_{k}}$ we place $s_{\beta_{i}}$ adjacent 
to $s_{\beta_{j}}$. The braid move 
$(s_{\beta_{i}},s_{\beta_{j}})\mapsto 
(s_{\beta_{i}}s_{\beta_{j}}s_{\beta_{i}},s_{\beta_{i}})=(s_{\lambda},s
_{\beta_{i}})$ yields a sequence with two reflections equal to 
$s_{\lambda}$. This contradicts Assumption~(*). If 
$\|\beta_{i}\|>\|\beta_{j}\|$ then 
$s_{\beta_{i}}(\beta_{j})=\beta_{j}+\beta_{i}$. This is a root with 
$(\lambda|\beta_{j}+\beta_{i})=1+\frac{1}{2}\|\beta_{i}\|\geq 3$. This 
is impossible by the choice of $\lambda$ in \S\ref{s2.14}. One reasons 
similarly in the case $\|\beta_{i}\|<\|\beta_{j}\|$ considering 
$s_{\beta_{j}}(\beta_{i})=\beta_{i}+\beta_{j}$. The claim of Step~2 is 
proved if $\beta_{1}=\lambda$. Let $\beta_{1}\neq \lambda$. Suppose 
$(\beta_{i}|\beta_{j})<0$ for some pair $i,j$ with $1\leq i<j\leq k$. 
The same arguments as above show that the only possibility is 
$\|\beta_{i}\|=\|\beta_{j}\|=2,\; \beta_{i}+\beta_{j}=\lambda$, in 
which case performing braid moves among the first $k$ reflections one 
obtains a Hurwitz system containing $s_{\lambda}$. One moves 
$s_{\lambda}$ to the first place by elementary transformations. The 
obtained Hurwitz system is of the type of Step~1 since none of 
$s_{\beta_{k+1}},\ldots,s_{\beta_{n}}$ has been changed with these 
transformations and by the minimality of $k$ (Condition~(iv) of 
Step~1). We already treated such cases in this step, so there is a 
braid-equivalent Hurwitz system for which Condition~(v) holds.

Step 3.\; We claim that for a Hurwitz system which satisfies 
conditions (i) -- (v) of Step~1 and Step~2 one of the following 
alternatives holds:

(vi)\; if $\beta_{1}\neq \lambda$ then $(\beta_{i}|\beta_{j})=0$ for 
$\forall i\neq j$ with $1\leq i,j\leq k$;

(vi)$'$\: if $\beta_{1}=\lambda$ then $(\beta_{i}|\beta_{j})=0$ for 
$\forall i\neq j$ with $2\leq i,j\leq k$. 

\noindent
Suppose $(\beta_{i}|\beta_{j})>0$ for some pair $i\neq j$ with $1\leq 
i,j\leq k$ and $\beta_{i}\neq \lambda\neq \beta_{j}$. We may assume 
$\|\beta_{i}\|\geq \|\beta_{j}\|$. If $\|\beta_{i}\|=\|\beta_{j}\|$ 
then $s_{\beta_{j}}(\beta_{i})=\beta_{i}-\beta_{j}$ and 
$(\lambda|s_{\beta_{j}}(\beta_{i}))=
\frac{1}{2}\|\beta_{i}\|-\frac{1}{2}\|\beta_
{j}\|=0$. If $\|\beta_{i}\|>\|\beta_{j}\|$ then $\|\beta_{j}\|=2$ and 
$(\lambda|\beta_{j})=1$ since $\lambda\neq \beta_{j}$. Furthermore 
$s_{\beta_{j}}(\beta_{i})=
\linebreak \beta_{i}-n(\beta_{i},\beta_{j})\beta_{j}$ 
and $(\lambda|s_{\beta_{j}}(\beta_{i})) = \frac{1}{2}\|\beta_{i}\|-
n(\beta_{i},\beta_{j}) = 0$. In both cases 
$(\lambda|s_{\beta_{j}}(\beta_{i}))=0$. Let $\beta$ be the positive 
root belonging to 
$\{s_{\beta_{j}}(\beta_{i}),-s_{\beta_{j}}(\beta_{i})\}$. Performing 
several braid moves among the first $k$ reflections we place 
$s_{\beta_{i}}$ adjacent to $s_{\beta_{j}}$. We obtain either 
$(\ldots,s_{\beta_{i}},s_{\beta_{j}},\ldots)\sim 
(\ldots,s_{\beta_{j}},s_{\beta},\ldots)$ or 
$(\ldots,s_{\beta_{j}},s_{\beta_{i}},\ldots)\sim 
(\ldots,s_{\beta},s_{\beta_{j}},\ldots)$ with $(\lambda|s_{\beta})=0$. 
In either case we move $s_{\beta}$ to the $k$-th place by successive 
braid moves. The obtained Hurwitz system contradicts the minimality of 
$k$ required in Condition~(iv) of Step~1.

Step 4.\; We claim Alternative~(vi) of Step~3 yields a contradiction. 
We have $\beta_{1},\ldots,\beta_{k}$ are $k$ mutually orthogonal 
roots. The involution $\sigma=s_{\beta_{1}}\cdots s_{\beta_{k}}$ has 
eigenvalues 1 and $-1$. Let $V_{+}$ and $V_{-}$ be the corresponding 
eigenspaces, $V=V_{+}\oplus V_{-}$. By \cite{Ca} Section~2 one has 
$V_{+}=H_{\beta_{1}}\cap \cdots \cap H_{\beta_{k}}$ and 
$V_{-}={V_{+}}^{\bot}$. The identity $s_{\beta_{1}}\cdots 
s_{\beta_{k}}s_{\beta_{k+1}}\cdots s_{\beta_{n}}=1$ and 
$(\lambda|\beta_{i})=0$ for $i\geq k+1$ implies 
$\sigma(\lambda)=\lambda$. Hence $\lambda \in V_{+}$, i.e. 
$(\lambda|\beta_{i})=0$ for $1\leq i\leq k$. This is a contradiction 
with Condition~(ii) of Step~1 (we recall $k\geq 1$ since 
$s_{\beta_{1}},\ldots,s_{\beta_{n}}$ generate $W$).

Step 5.\; If $R$ is not of type $G_{2}$ then Alternative~(vi)$'$ of 
Step~3 yields a contradiction. First we prove $k\leq 2$ in this case. 
Suppose among $\beta_{2},\ldots,\beta_{k}$ there are two roots with 
length 2. We may assume without loss of generality that 
$\|\beta_{2}\|=\|\beta_{3}\|=2$. We have 
$s_{\beta_{2}}(\lambda)=\lambda-\beta_{2}$ and 
$s_{\beta_{3}}(s_{\beta_{2}}(\lambda))=\lambda-\beta_{2}-\beta_{3}=
\beta$. The latter root satisfies $(\lambda|\beta)=0$. Performing the 
braid moves $(s_{\lambda},s_{\beta_{1}},s_{\beta_{2}},\ldots)\sim 
(s_{\beta_{1},s_{\lambda-\beta_{1}}},s_{\beta_{2}},\ldots )\sim 
(s_{\beta_{1}},s_{\beta_{2}},s_{\beta},\ldots )$ and moving 
$s_{\beta}$ to the $k$-th place by elementary transformations we 
obtain a Hurwitz system with $n-k+1$ reflections which fix $\lambda$. 
This contradicts Condition~(iv) of Step~1. Therefore among 
$\beta_{2},\ldots,\beta_{k}$ there is at most one root of length 2. If 
$R$ is of type $B_r,C_r$ or $F_4$ suppose among 
$\beta_{2},\ldots,\beta_{k}$ there is a long root. We may assume 
without loss of generality this is $\beta_{2}$. Then 
$s_{\beta_{2}}(\lambda)=\lambda-\beta_{2}=\beta$ and 
$(\lambda|\beta)=(\lambda|\lambda)-(\lambda|\beta_{2})=2-\frac{1}{2}
\|\beta_{2}\|=0$. Performing the braid move 
$(s_{\lambda},s_{\beta_{2}},\ldots)\sim 
(s_{\beta_{2}},s_{\beta},\ldots)$ and moving $s_{\beta}$ to the $k$-th 
place by successive elementary transformations we obtain again a 
Hurwitz system which contradicts the minimality of $k$. Our claim that 
$k\leq 2$ is proved. We have 
$s_{\lambda}=s_{\beta_{2}}s_{\beta_{3}}\cdots s_{\beta_{n}}$.
If $k=1$ then $(\lambda|\beta_{i})=0$ for $i\geq 2$. Applying
both sides of the above equality to $\lambda$ we obtain the absurdity 
$-\lambda =\lambda$.
 If $k=2$ 
then  
$(\lambda|\beta_{i})=0$ for $i\geq 3$. Applying both sides of the 
above equality to 
$\lambda$ we obtain 
$-\lambda=s_{\lambda}(\lambda)=s_{\beta_{2}}(\lambda)$. Therefore 
$\beta_{2}=\lambda$ which contradicts Assumption~(*).

Step 6.\; It remains to consider Alternative~(vi)$'$ of Step~3 when 
$R$ is of type $G_{2}$. The same argument as in Step~5 shows that 
among $\beta_{2},\ldots,\beta_{k}$ there is at most one root of length 
2. We claim there is also at most one long root 
among $\beta_{2},\ldots,\beta_{k}$. Indeed if $\beta_{i},\beta_{j}$  
are long roots, then $s_{\beta_{i}}(\lambda)=\lambda-\beta_{i},\; 
s_{\beta_{j}}(s_{\beta_{i}}(\lambda))=\lambda-\beta_{i}-\beta_{j}$. 
The root $\lambda-\beta_{i}-\beta_{j}$ satisfies 
$(\lambda|\lambda-\beta_{i}-\beta_{j})=2-3-3=-4$. This is an absurdity 
since for $\forall \gamma \in R$ one has $(\lambda|\gamma)\in 
\{0,\pm1,\pm3\}$. The case $k\leq 2$ is impossible by the argument of 
Step~5. It remains to consider the case where $k=3$ and 
$\{\beta_{2},\beta_{3}\}$ consists of a short and a long root which 
are orthogonal. Without loss of generality we may assume 
$\|\beta_{2}\|=2,\; \|\beta_{3}\|=6$. The number $n$ is even, so 
$n\geq 4$. If among $\beta_{4},\ldots,\beta_{n}$ there were no long 
roots, then $s_{\beta_{3}}$ would be contained in the subgroup 
generated by reflections with respect to short roots. This is 
impossible since the only reflections in the latter subgroup  are 
$s_{\gamma}$ with $\gamma$ a short root in $R$. Without loss of 
generality we may assume $\beta_{4}$ is a long root. We obtain a 
Hurwitz system 
$(s_{\lambda},s_{\beta_{2}},s_{\beta_{3}},s_{\beta_{4}},\ldots)$ 
braid-equivalent to $t_{1}\cdots t_{n}=1$ where 
$\|\beta_{2}\|=2, \|\beta_{3}\|=\|\beta_{4}\|=6, 
(\lambda|\beta_{2})=1, (\lambda|\beta_{3})=3, (\lambda|\beta_{4})=0, 
(\beta_{2}|\beta_{3})=0$. Table X of \cite{Bo} Ch.VI lists the 
positive roots in $R$ of type $G_{2}$: \; $R^{+}=R^{+}_{s}\cup 
R^{-}_{_{\ell}}, \; R^{+}_{s}=\{\alpha_{1},\alpha'_{1},\omega_{1}\},\; 
R^{-}_{_{\ell}} = \{\alpha_{2},\alpha'_{2},\omega_{2}\}$ where 
$\alpha_{1},\alpha_{2}$ is a base of $R$ and $\omega_{1},\omega_{2}$ are the 
fundamental weights. We have $\lambda=\omega_{1}$. The following two 
possibilities for 
$(s_{\lambda},s_{\beta_{2}},s_{\beta_{3}},s_{\beta_{4}})$ may occur:
\begin{equation}\label{es2.22}
 (s_{\omega_{1}},s_{\alpha_{1}},s_{\omega_{2}},s_{\alpha_{2}})\quad \quad 
\text{or}\quad \quad
(s_{\omega_{1}},s_{\alpha'_{1}},s_{\alpha'_{2}},s_{\alpha_{2}}).
\end{equation}
In the first case we perform the following braid moves:
\begin{equation*}
\begin{split}
&(s_{\omega_{1}},s_{\alpha_{1}},s_{\omega_{2}},s_{\alpha_{2}})\sim
(s_{\omega_{1}},s_{\omega_{2}},s_{\alpha_{1}},s_{\alpha_{2}})\sim
(s_{\omega_{1}},s_{\omega_{2}},s_{\alpha_{2}},s_{\alpha'_{1}})\\
\sim
&(s_{\omega_{1}},s_{\omega_{2}},s_{\alpha'_{1}},s_{-\omega_{2}})
\sim
(s_{\omega_{1}},s_{-\omega_{1}},s_{\omega_{2}},s_{-\omega_{2}})=
(s_{\omega_{1}},s_{\omega_{1}},s_{\omega_{2}},s_{\omega_{2}}).
\end{split}
\end{equation*}
This contradicts Assumption~(*). The product of the first quadruple of 
\eqref{es2.22} 
equals 1 and conjugating it by $s_{\omega_{1}}$ one obtains the second
quadruple. By Lemma~\ref{s1.9aa} the second quadruple is 
braid-equivalent to the first one, hence braid-equivalent to 
$(s_{\omega_{1}},s_{\omega_{1}},s_{\omega_{2}},s_{\omega_{2}})$. 
This contradicts Assumption~(*).

We proved that Assumption~(*) leads to a contradiction. Therefore the 
Hurwitz system $t_{1}\cdots t_{n}=1$ is braid-equivalent to some 
$t'_{1}t'_{2}\cdots t'_{n}$ with $t'_{1}=t'_{2}$. 
The proposition is proved.
\end{proof}
\begin{pro}\label{s2.23}
Let $R$ be a root system and let $W$ be its Weyl group. Let 
$t_{1}\cdots t_{n}=1$ be a Hurwitz system of reflections in $W$. Then 
$t_{1}\cdots t_{n}=1$ is braid-equivalent to a Hurwitz system 
$t'_{1}\cdots t'_{n}=1$ where $t'_{2i-1}=t'_{2i}$ for 
$i=1,\ldots,\frac{n}{2}$. 
\end{pro}
\begin{proof}
Use Proposition~\ref{s2.15} and induction on the even number $n$.
\end{proof}
\begin{block}\label{s2.24}
Let $R$ be an irreducible root system of rank $r$ with Weyl group $W$. 
In \S\ref{s1.7a} we defined the Hurwitz spaces 
$\mathcal{H}_{W;\underline{n}}(\mathbb{P}^1,b_0)$ and 
$\mathcal{H}_{W;\underline{n}}(\mathbb{P}^1)$. If $R$ is simply laced, 
i.e. of type $A_r,D_r,E_6,E_7$ or $E_8$, one has $\underline{n}=n$ and 
the spaces parametrize irreducible Galois covers of $\mathbb{P}^1$ 
branched in $n$ points whose local monodromies are reflections in $W$. 
If $R$ is non-simply laced, i.e. of type $B_r,C_r,F_4$ or $G_2$, one 
has $\underline{n}=(n_{s},n_{_{\ell}})$ and the spaces parametrize 
irreducible Galois covers of $\mathbb{P}^1$ branched in 
$n=n_{s}+n_{_{\ell}}$ points with $n_{s}$ discriminant points whose local 
monodromies are reflections with respect to short roots and 
$n_{_{\ell}}$ discriminant points whose local 
monodromies are reflections with respect to long roots. In the 
non-simply laced case let $r_{s}$, respectively $r_{_{\ell}}$, denote the 
number of short roots, respectively long roots, in the Dynkin diagram 
of $R$. One has $r=r_{s}+r_{_{\ell}}$ and for types $B_r,C_r,F_4$ and
$G_2$ the pair $(r_{s},r_{_{\ell}})$ equals respectively $(1,r-1),\: 
(r-1,1),\: (2,2)$ and $(1,1)$ (cf. \cite{Bo} Ch.VI).
\end{block}
\begin{thm}\label{s2.25}
Let $R$ be an irreducible root system of rank $r$  with Weyl group 
$W$. 
\renewcommand{\theenumi}{\roman{enumi}}
\begin{enumerate}
\item 
The Hurwitz spaces 
$\mathcal{H}_{W;\underline{n}}(\mathbb{P}^1,b_0)$ and 
$\mathcal{H}_{W;\underline{n}}(\mathbb{P}^1)$ are irreducible when 
non-empty.
\item 
If $R$ is of type $A_r,D_r,E_6,E_7$ or $E_8$ then 
$\mathcal{H}_{W;\underline{n}}(\mathbb{P}^1,b_0)\neq \emptyset$, 
equivalently $\mathcal{H}_{W;\underline{n}}(\mathbb{P}^1)\neq 
\emptyset$, if and only if $n\geq 2r$.
\item 
If $R$ is of type $B_r,C_r,F_4$ or $G_2$ then 
$\mathcal{H}_{W;\underline{n}}(\mathbb{P}^1,b_0)\neq \emptyset$, 
equivalently $\mathcal{H}_{W;\underline{n}}(\mathbb{P}^1)\neq 
\emptyset$, if and only if $n_{s}\equiv 0(mod\, 2),\; n_{_{\ell}}\equiv 
0(mod\, 2),\; n_{s}\geq 2r_{s}$ and $n_{_{\ell}}\geq 2r_{_{\ell}}$.
\end{enumerate}
\end{thm}
\begin{proof}
By \S\ref{s1.7a}(i) and Corollary~\ref{s1.9ad} it suffices to prove 
the statements for $\mathcal{H}_{W;\underline{n}}(\mathbb{P}^1,b_0)$. 
Let $\alpha_{1},\ldots,\alpha_{r}$ be a base of $R$. If $R$ is simply 
laced let us choose an arbitrary root $\alpha$. If $R$ is 
non-simply laced let us choose an arbitrary short root $\alpha$ and an
arbitrary long root $\beta$. According to \S\ref{s1.9ac} the following 
claim proves Part~(i) and one of the directions of Part~(ii) and 
Part~(iii).

CLAIM. 
\emph{Let $t_{1}\cdots t_{n}=1$ be a Hurwitz system of reflections in 
$W$ such that $t_1,\ldots,t_n$ generate $W$. If $R$ is simply laced 
then $n\geq 2r$ and the Hurwitz system is braid-equivalent to 
\begin{equation}\label{es2.26a}
s_{\alpha_{1}}s_{\alpha_{1}}s_{\alpha_{2}}s_{\alpha_{2}}\cdots 
s_{\alpha_{r}}s_{\alpha_{r}}s_{\alpha}\cdots s_{\alpha}=1
\end{equation}
where $s_{\alpha_{}}$ appears $n-2r$ times. If $R$ is non-simply laced 
then $n_{s}$ and $n_{_{\ell}}$ are even, $n_{s}\geq 2r_{s},n_{_{\ell}}\geq 
2r_{_{\ell}}$ and the Hurwitz system is braid-equivalent to 
\begin{equation}\label{es2.26b}
s_{\alpha_{1}}s_{\alpha_{1}}s_{\alpha_{2}}s_{\alpha_{2}}\cdots 
s_{\alpha_{r}}s_{\alpha_{r}}s_{\alpha}\cdots s_{\alpha}s_{\beta}\cdots 
s_{\beta} = 1
\end{equation}
where $s_{\alpha}$ appears $n_{s}-2r_{s}$ times and $s_{\beta}$ 
appears $n_{_{\ell}}-2r_{_{\ell}}$ times.}

Let $n=2m$. According to Proposition~\ref{s2.23} the Hurwitz system 
$t_{1}\cdots t_{n}=1$ is braid-equivalent to $t'_{1}\cdots t'_{n}=1$ 
where $t'_{2i-1}=t'_{2i}=s_{\beta_{i}}$ for $i=1,\ldots,m$. One has 
$\langle s_{\beta_{1}},\ldots,s_{\beta_{m}}\rangle = \langle 
t_{1},\ldots,t_{n}\rangle=W$. According to Corollary~\ref{s1.13a} 
there is a finite sequence of Nielsen transformations by which one can 
obtain $\{s_{\alpha_{1}},\ldots,s_{\alpha_{r}}\}$ from 
$\{s_{\beta_{1}},\ldots,s_{\beta_{m}}\}$. Lemma~\ref{s1.9b} shows that 
if a Hurwitz system contains two adjacent pairs of involutions $(s,s)$ 
and $(t,t)$ then replacing $(t,t)$ by $(sts,sts)$ one obtains a 
braid-equivalent Hurwitz system. This implies that extending the above 
sequence  of Nielsen transformations to Nielsen transformations of 
Hurwitz systems composed of pairs of elements of the corresponding 
sets one obtains a sequence of braid-equivalences. Eventually we 
obtain a Hurwitz system composed of pairs 
$(s_{\alpha_{i}},s_{\alpha_{i}})$ and every such pair with $1\leq i\leq 
r$ do appear. Using Lemma~\ref{s1.9ab} we may replace the obtained 
Hurwitz system by a braid-equivalent one in which the first $2r$ 
elements are 
$(s_{\alpha_{1}},s_{\alpha_{1}},
\ldots,s_{\alpha_{r}},s_{\alpha_{r}})$. These $2r$ reflections 
generate $W$, so by Lemma~\ref{s1.9b} we may replace any of the 
remaining pairs $(s_{\gamma},s_{\gamma})$ by $(s_{\alpha},s_{\alpha})$ 
if $\|\gamma\|=\|\alpha\|$ and by $(s_{\beta},s_{\beta})$ if 
$\|\gamma\|=\|\beta\|$. This proves the claim.

For the proof of the other direction of Part~(ii) and Part~(iii) 
notice that  if the specified inequalities for the number of 
discriminant points are valid one may define Hurwitz systems by 
\eqref{es2.26a} and \eqref{es2.26b} and apply \S\ref{s1.7a}(i).
\end{proof}
Performing a sequence of braid moves to a given Hurwitz system may be 
viewed in two ways. Either one fixes a simple system of closed arcs in 
$\mathbb{P}^1-D$ based at $b_{0}$ and varies the homomorphism 
$m:\pi_{1}(\mathbb{P}^1-D,b_{0})\to W$ thus obtaining information 
about the connected components of the Hurwitz space 
$\mathcal{H}_{W;\underline{n}}(\mathbb{P}^1,b_0)$, or one fixes the 
monodromy map \linebreak
$m:\pi_{1}(\mathbb{P}^1-D,b_{0})\to W$ and varies the simple arc 
system thus obtaining a normal form for the local monodromies of a 
given covering and eventually determining the topological type of the 
covering. So the proof of the theorem yields the following result in 
which we use the notation introduced in \S\ref{s2.24}.
\begin{cor}\label{s2.29}
Let $\pi : X\to \mathbb{P}^1$ be an irreducible Galois cover with 
Galois group the Weyl group of an irreducible root system $R$ of rank 
$r$. Suppose $\pi$ is simply ramified, i.e. every local monodromy is a 
reflection. Let $D$ be the discriminant locus of $\pi$, let $|D|=n$, 
let $b_{0}\in \mathbb{P}^1-D$, let $x_{0}\in \pi^{-1}(b_{0})$ and let 
$m=m_{x_{0}}:\pi_{1}(\mathbb{P}^1-D,b_{0})\to W$ be the monodromy map. 
Then there is a simple arc system with initial point $b_{0}$ and end 
points in $D$ such that the local monodromies along the corresponding 
simple system of closed arcs $\gamma_{1},\ldots,\gamma_{n}$ have the  
form given by the Hurwitz system \eqref{es2.26a}, respectively 
\eqref{es2.26b}. Namely fixing a basis $\alpha_{1},\ldots,\alpha_{r}$ 
of $R$, an arbitrary root $\alpha$ if $R$ is simply laced, an 
arbitrary short root $\alpha$ and an arbitrary long root $\beta$ if 
$R$ is non-simply laced, one  has:
\renewcommand{\theenumi}{\roman{enumi}}
\begin{enumerate}
\item if $R$ is of type $A_r,D_r,E_6,E_7$ or $E_8$ then 
\[
(m(\gamma_{1}),\ldots,m(\gamma_{n})) = 
(s_{\alpha_{1}},s_{\alpha_{1}},s_{\alpha_{2}},s_{\alpha_{2}},\cdots, 
s_{\alpha_{r}},s_{\alpha_{r}},s_{\alpha},\cdots,s_{\alpha})
\]
where $s_{\alpha}$ appears $n-2r$ times;
\item if $R$ is of type $B_r,C_r,F_4$ or $G_2$ then 
\[
(m(\gamma_{1}),\ldots,m(\gamma_{n})) = 
(s_{\alpha_{1}},s_{\alpha_{1}},s_{\alpha_{2}},s_{\alpha_{2}},\cdots, 
s_{\alpha_{r}},s_{\alpha_{r}},s_{\alpha},\cdots,s_{\alpha},s_{\beta},\cdots, 
s_{\beta})
\]
where $s_{\alpha}$ appears $n_{s}-2r_{s}$ times and $s_{\beta}$ 
appears $n_{_{\ell}}-2r_{_{\ell}}$ times.
\end{enumerate}
\end{cor}
We now extend Theorem~\ref{s2.25} and Corollary~\ref{s2.29} to simply 
ramified $W$-coverings where $W$ is an arbitrary Weyl group. We refer 
to \S\ref{s1.9da} and \S\ref{s2.24} for the notation used. The 
superscript $^{(i)}$ refers to the irreducible root system $R^{(i)}$. 
\begin{thm}\label{s2.31}
Let $R$ be a root system  with Weyl group 
$W$. Let $R=R^{(1)}\sqcup \ldots \sqcup R^{(k)}$ be its decomposition 
into irreducible components.
\renewcommand{\theenumi}{\roman{enumi}}
\begin{enumerate}
\item 
The Hurwitz spaces 
$\mathcal{H}_{W;\underline{n}}(\mathbb{P}^1,b_0)$ and 
$\mathcal{H}_{W;\underline{n}}(\mathbb{P}^1)$ are irreducible when 
non-empty.
\item 
$\mathcal{H}_{W;\underline{n}}(\mathbb{P}^1,b_0)\neq \emptyset$, 
equivalently $\mathcal{H}_{W;\underline{n}}(\mathbb{P}^1)\neq 
\emptyset$, if and only if 
for every simply laced component $R^{(i)}$ one has $n^{(i)}\equiv 0
(mod\, 2),\; n^{(i)}\geq 2r^{(i)}$ and for every non-simply laced 
component $R^{(j)}$ one has $n_{s}^{(j)}\equiv 0(mod\, 2),\; 
n_{_{\ell}}^{(j)}\equiv 0(mod\, 2),\; n_{s}^{(j)}\geq 2r_{s}^{(j)},\; 
n_{_{\ell}}^{(j)}\geq 2r_{_{\ell}}^{(j)}$.
\end{enumerate}
\end{thm}
\begin{proof}
Let $t_{1}\cdots t_{n}=1$ be a Hurwitz system of reflections in $W$ 
such that $t_1,\ldots,t_n$ generate $W$. It is braid-equivalent to a 
concatenation $T^{(1)}T^{(2)}\cdots T^{(k)}=1$ where each $T^{(i)}$ is 
a Hurwitz system of reflections in $W^{(i)}$ which generate $W^{(i)}$ 
(cf. the proof of Proposition~\ref{s2.15}). The theorem is proved 
applying the argument of Theorem~\ref{s2.25} to each $T^{(i)}$.
\end{proof}
\begin{cor}\label{s2.32}
Let $\pi : X\to \mathbb{P}^1$ be an irreducible Galois cover with 
Galois group the Weyl group of an arbitrary root system $R$. Suppose 
every local monodromy is a reflection. Let $D$ be the discriminant 
locus of $\pi$, let $b_{0}\in \mathbb{P}^1-D$ and let $D=
D^{(1)}\sqcup \ldots \sqcup D^{(k)}$ be the disjoint union 
corresponding to the decomposition into irreducible components $R=
R^{(1)}\sqcup \ldots \sqcup R^{(k)}$. Let $n^{(i)}=|D^{(i)}|$. Then 
there is a simple arc system with initial point $b_{0}$ and end points 
in $D$, ordered so that the first $n^{(1)}$ arcs end in $D^{(1)}$, the 
arcs with numbers $n^{(1)}+1,\ldots,n^{(1)}+n^{(2)}$ end in $D^{(2)}$ 
etc., such that for every $i$ the local monodromies corresponding to 
the collection of arcs ending in $D^{(i)}$ are given by the formulae in 
Corollary~\ref{s2.29} with $R$ replaced by $R^{(i)}$.
\end{cor}

\bibliographystyle{amsalpha}
\providecommand{\bysame}{\leavevmode\hbox to3em{\hrulefill}\thinspace}

\bigskip

\noindent
{\sc
Dipartimento di Matematica ed Applicazioni\\
Universit\`{a} degli Studi di Palermo\\
Via Archirafi n. 34, 90123 Palermo, Italy\\
and \\
Institute of Mathematics, Bulgarian Academy of Sciences}

\medskip \noindent
{\it E-mail address:} {\bf kanev@math.unipa.it}

\end{document}